\documentclass[a4paper,12pt]{article}

\usepackage[latin1]{inputenc}
\usepackage[english]{babel}
\usepackage{amssymb}
\usepackage{amsmath}
\usepackage{latexsym}
\usepackage{amsthm}
\usepackage[pdftex]{graphicx}
\usepackage{tikz} %questo package va messo dopo graphicx
\usepackage{pgfplots}  
\usepackage{mathtools}
\pgfplotsset{width=6.6cm,compat=1.7}
\usepackage{tkz-euclide}

\usetikzlibrary{cd}

\usetikzlibrary{decorations.pathreplacing}

 \theoremstyle{plain}
 \newtheorem{thm}{Theorem}[section]
 \newtheorem{cor}[thm]{Corollary}
  \newtheorem{conj}[thm]{Conjecture}
 \newtheorem{lem}[thm]{Lemma}
 \newtheorem{prop}[thm]{Proposition}

 \theoremstyle{definition}
 
 \newtheorem{example}[thm]{Example}

 \theoremstyle{remark}
 
\title{On the number of fixed points of the map $\gamma$}
\date{}
\author{
Niccol\`o Castronuovo \\
\texttt{castronuovoniccolo@gmail.com}
}
\begin{document}
\maketitle

\begin{abstract}
We recursively define a sequence $\{F_{n,k}\}_{n,k\in\mathbb N }$ and we prove that such sequence contains only the symbols $\{0,1\}.$
We investigate some number-theoretic properties of such sequence and of the way it can be generated.
The number $F_n$ can be interpreted as the number of  fixed points of semilength $n$ of the map $\gamma$ introduced in \cite{BaBoCaCo}. Our results partially answer conjectures posed to the author by Cori \cite{test1}.
\end{abstract}

%\noindent {\bf Keywords:}  

%\noindent {\bf MSC2010:} 05A05, 05A15 (primary); 05A10 (secondary).

\section{Introduction}

In this paper we consider an infinite $(0,1)$-matrix $F$ defined  in the following way. Let $F:=[F_{n,k}]_{n\geq 0\,k\geq 0}$ be the  doubly-infinite matrix all of whose entries are equal to $0.$
Apply to $F$ the following step:
\begin{itemize}
    \item[Step 0] Set $F_{0,0}=1.$
\end{itemize}
For all $i\geq 1$ apply to $F$ the following step:
\begin{itemize}
    \item [Step $i$] For each pair $(n,k)$ such that the entry $F_{n,k}$ changes its value in Step $i-1,$ increase $F_{n+k,k}$
    and $F_{3n+1-2k,2n+1-k}$ by $1.$
\end{itemize}

$F$ is the matrix obtained in this way.

We say that an entry $F_{n,k}$ of matrix $F$ \textit{is created} at Step $i$ if $F_{n,k}>0$ and, during the creation of matrix $F$ it changes its value during Step $i.$

It is trivial to verify that $F_{0,k}>0$ if and only if $k=0$ and that $F_{n,k}=0$ if $k>n.$ Hence the matrix $F$ is lower triangular and $\{F_{n,k}\}_{n\geq 0\,k\geq 0}$ can be thought as a doubly-indexed sequence.

The matrix $F$ is related to the map $\gamma,$
a bijection defined over the set of Dyck words of semilength  $n.$ This map and its properties are introduced in  \cite{BaBoCaCo} and further studied in \cite{CasCo} and \cite{test1}. This last paper, in particular, deals with the characterization of the fixed points of $\gamma.$
%Here we retain the same definitions of these papers.

More precisely, $F_{n,k}$ is equal to the number of Dyck words of semilength $n,$ with  principal prefix of length $k$ and fixed under the action of $\gamma.$
The fact that there is at most one of such words (see \cite{BaBoCaCo}) implies that the matrix $F$ is a $0-1$ matrix. We will reprove this result in Corollary \ref{zeroone}.
The sum of entries in row $n,$ $F_n:=\sum_k F_{n,k},$  is the total number of Dyck words of semilength $n$ fixed by $\gamma$ (see \cite{BaBoCaCo} for the main definitions).

The reason for which we do not reintroduce the definition of the map $\gamma$ is that the sequence $\{F_{n,k}\}_{n,k}$ can be defined implicitly as above (see \cite{test1}). Hence all the results of the paper can be stated in a number-theoretic form without appealing to the original definition of $F_{n}.$

The first few rows of the matrix are reported below (the elements above the main diagonal are all zeros and are not indicated).

\setcounter{MaxMatrixCols}{20}

$$
F=
 \begin{matrix}
   1 &   &   &   &   &   &   &   &   &   &   &   &   &   &   &   &  &   \\
   0 & 1 &   &   &   &   &   &   &   &   &   &   &   &   &   &   &   &  \\
   0 & 1 & 1 &   &   &   &   &   &   &   &   &   &   &   &   &   &   &  \\
   0 & 1 & 0 & 1 &   &   &   &   &   &   &   &   &   &   &   &   &   &  \\
   0 & 1 & 1 & 0 & 1 &   &   &   &   &   &   &   &   &   &   &   &   &  \\
   0 & 1 & 0 & 0 & 1 & 1 &   &   &   &   &   &   &   &   &   &   &   &  \\
   0 & 1 & 1 & 1 & 0 & 0 & 1 &   &   &   &   &   &   &   &   &   &   &  \\
   0 & 1 & 0 & 0 & 0 & 0 & 0 & 1 &   &   &   &   &   &   &   &   &   &  \\
   0 & 1 & 1 & 0 & 1 & 0 & 1 & 1 & 1 &   &   &   &   &   &   &   &   &  \\
   0 & 1 & 0 & 1 & 1 & 0 & 0 & 1 & 0 & 1 &   &   &   &   &   &   &   &  \\
   0 & 1 & 1 & 0 & 0 & 1 & 0 & 0 & 0 & 0 & 1 &   &   &   &   &   &   &  \\
   0 & 1 & 0 & 0 & 0 & 0 & 0 & 0 & 1 & 0 & 1 & 1 &   &   &   &   &   &  \\
   0 & 1 & 1 & 1 & 1 & 0 & 1 & 0 & 0 & 0 & 0 & 0 & 1 &   &   &   &   &  \\
   0 & 1 & 0 & 0 & 1 & 0 & 0 & 0 & 0 & 0 & 1 & 1 & 0 & 1 &   &   &   &  \\ 
   0 & 1 & 1 & 0 & 0 & 0 & 1 & 1 & 0 & 0 & 1 & 0 & 1 & 1 & 1 &   &   &  \\ 
   0 & 1 & 0 & 1 & 0 & 1 & 0 & 1 & 0 & 0 & 0 & 1 & 0 & 0 & 0 & 1 &   &  \\  
   0 & 1 & 1 & 0 & 1 & 0 & 0 & 1 & 1 & 0 & 0 & 0 & 0 & 0 & 0 & 0 & 1 &  \\ 
   0 & 1 & 0 & 0 & 1 & 0 & 0 & 0 & 0 & 0 & 0 & 0 & 1 & 1 & 0 & 0 & 1 & 1 \\
   & \dots 
   \end{matrix}
$$

The first values of the sequence $\{F_n\}_{n\geq 0}$ are
$$1,1,2,2,3,3,4,2,6,5,4,4,6,5,8,6,6,6\ldots$$

We study the matrix $F$ and the sequence $F_n,$ and investigate their properties.

Our results partially answer the following two conjectures posed to the author by Cori \cite{test1}.
\begin{conj}\label{conj1}
$F_n\geq 3$ for all $n> 7.$
\end{conj}
\begin{conj}\label{conj2}
$\lim_{n\to \infty} F_n=\infty.$
\end{conj}

In particular we answer in the affirmative Conjecture \ref{conj1} and give some results toward the solution to Conjecture \ref{conj2}.

\section{The matrix $F$ and a free subsemigroup of $SL(3,\mathbb Z).$}

Identify the entry $F_{i,j}$ of the matrix $F$ with the integer vector with coordinates $(i,j)$ in the $\mathbb Z\times \mathbb Z$ lattice plane. 
It follows immediately  from the definition of the matrix $F,$ that, if $i,j>1,$  $F_{i,j}>0$ if and only if the vector $(i,j)$ can be reached iteratively applying to the vector $(1,1)$ the following affine transformations (in arbitrary order)
$$
\hat G:\begin{pmatrix}x\\ y
\end{pmatrix} \to \begin{pmatrix}1 & 1\\ 0 & 1
\end{pmatrix}\cdot \begin{pmatrix}x\\ y
\end{pmatrix},
$$
and
$$
\hat S:\begin{pmatrix}x\\ y
\end{pmatrix} \to \begin{pmatrix}3 & -2\\2 & -1
\end{pmatrix}\cdot \begin{pmatrix}x\\ y
\end{pmatrix}+\begin{pmatrix}1\\ 1
\end{pmatrix}.
$$

The group generated by the affine transformations $\hat G$ and $\hat S$ can be identified with the subgroup $\langle S,G\rangle$ of $SL(3,\mathbb Z)$ generated by matrices
$$G= \begin{pmatrix}1 & 1 & 0 \\ 0 & 1 & 0 \\ 0 & 0 & 1
\end{pmatrix}\quad \mbox{ and }
\quad S= \begin{pmatrix} 3 & -2 & 1 \\ 2 &-1 &1 \\ 0 & 0 & 1
\end{pmatrix}. $$
In this identification, the actions of $\hat S$ and $\hat G$ on the lattice point $(x,y)\in \mathbb Z^2$ correspond to the actions of $S$ and $G,$ respectively, on the lattice point $(x,y,1) \in \mathbb Z^3.$
Hence we can consider directly the action of the group $<S,G>$ over the set $\mathbb Z^2.$ If $w\in <S,G>,$ we denote by $w(x,y)$ the image of the vector $(x,y)$ under this action. 

Now we study some properties of the group $<S,G>$ and of matrices $S$ and $G.$

Notice that the group generated by $S$ and $G$ is not free. In fact we have $(GS^{-1})^4=\begin{pmatrix}1 & 0 & 0 \\ 0 & 1 & 0 \\ 0 & 0 & 1
\end{pmatrix}.$ If we restrict our attention to the monoid generated by $S$ and $G,$ it is free.

\begin{thm}\label{free}
The monoid $H$ generated by $S$ and $G$ is free. The action of $H$ on $\mathbb Z^2$ is free i.e. for all $(x,y)\in \mathbb Z^2$ if $w(x,y)=w'(x,y),$ with $w,w'\in H,$ then $w=w'.$
\end{thm}

For the proof it will be useful a version of the so called \textit{Ping-Pong lemma} or \textit{Table-Tennis lemma} for semigroups (see e.g. \cite{de2000topics}[p. 188]), that we report here.

\begin{lem}
Let $\Gamma$ be a group acting on a set $X.$
Assume that there exist $\gamma_1,\gamma_2\in \Gamma$ and $X_1,X_2\subseteq X$ such that $X_1\cap X_2=\emptyset,$
$\gamma_1(X_1\cup X_2)\subseteq X_1$ and $\gamma_2(X_1\cup X_2)\subseteq X_2.$ Then the semigroup generated in $\Gamma$ by $\gamma_1$ and $\gamma_2$ is free.
\end{lem}

Now we proceed to the proof of Theorem \ref{free}.
\proof
For the first part, we apply the previous lemma with $\gamma_1:=G,$ $\gamma_2:=S$ and $\Gamma:=H.$
As described above, $H$ acts in the standard way on the set $X=\mathbb Z^2,$ and, more generally, on $\mathbb R^2.$
Notice that the only fixed point under the action of $H$ on $\mathbb R^2$ is $(-0.5,0),$ hence any line through this point is mapped onto another such line. Moreover, every point of the line $y=x+\frac{1}{2}$ is fixed by $S$ and every point of the line $y=0$ is fixed by $G.$

We consider the following disjoint subsets of $X:$
$$X_1:=\{(x,y)\in \mathbb Z^2\,|\, x>-\frac{1}{2}\,,\,0<y<\frac{x}{2}+\frac{1}{4}\} $$
$$X_2:=\{(x,y)\,|\, x>-\frac{1}{2}\,,\,\frac{2x}{3}+\frac{1}{3}<y<x+\frac{1}{2}\} $$

Those subsets are depicted in the figure below.

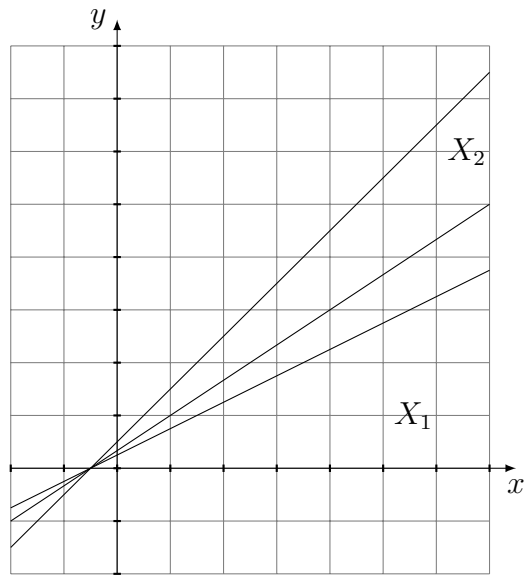
\begin{center}
\begin{tikzpicture}[scale=0.7] 
  \tkzInit[xmax=7,ymax=8,xmin=-2,ymin=-2]
   \tkzGrid
\draw[->,ultra thick] (-2,0)--(8,0) node[right]{$x$};
\draw[->,ultra thick] (0,-2)--(0,9) node[above]{$y$};
   % \tkzDrawXY
    %yticklabels={,,}
 \draw[color=black,domain=-2:7]    plot (\x,\x+1/2)             node[right] {}; 
    \draw[color=black,domain=-2:7] plot (\x,\x/2+1/4) node[right] {};
    \draw[color=black,domain=-2:7] plot (\x,2*\x/3+1/3) node[right] {};    
    \draw (5,1) node[right] {$X_1$};
    \draw (6,6) node[right] {$X_2$};    
  \end{tikzpicture}
\end{center}

It is trivial to verify that 
$$G(X_1\cup X_2)\subseteq X_1,\quad\quad S(X_1\cup X_2)\subseteq X_2.$$

Now we prove the second part of the theorem.
We want to show that $$\forall (x,y)\in \mathbb Z^2,\quad \mbox{ if }w(x,y)=w'(x,y),\mbox{ with }w,w'\in H,\,\mbox{ then }w=w'.$$
We proceed by induction on the minimum of the lengths of $w$ and $w',$ thought as words in the letters $S$ and $G,$ $$m:=\min \{|w|,|w'|\}.$$
If $m=0$ then one of the two words is the identity. The equation $w(x,y)=(x,y),$ with $w$ different from the identity is clearly impossible
 since both $S$ and $G$ increase the abscissa of the point on which they act. 
 Suppose the assertion true for all values of $m$ up to $N.$
If $m=N+1$ and the first letter of $w$ and $w'$ is the same, e.g. the letter $G,$ we have $w=G \hat w,$ $w'=G \hat w'$
and $G \hat w (x,y)=G \hat w'(x,y).$ This implies $\hat w (x,y)=\hat w'(x,y)$  and hence $w=w'$ by the inductive hypothesis. 
If $m>0$ and $w=G\hat w$ and $w'=S\hat{w'},$
then $w(x,y)\in X_1$ and $w'(x,y)\in X_2.$ Since $X_1\cap X_2=\emptyset$ it is impossible that $w(x,y)=w'(x,y).$

This concludes the proof. 
\endproof

As recalled above, there is a bijection between $H$ and the set $\{(i,j)|\,i,j>0,\, F_{i,j}\neq 0\}.$ This bijection maps an element $w\in H$ to the pair $w(1,1).$
Hence the previous theorem leads to the following corollary.

\begin{cor}\label{zeroone}
The matrix $F$ is a 0-1 matrix. 
\end{cor}

The following lemma will be useful in the sequel.

\begin{lem}\label{matrices}
The matrices $S$ and $G$ satisfy the following identities for all $i$
$$G^i=\begin{pmatrix}1 & i & 0 \\ 0 & 1 & 0 \\ 0 & 0 & 1
\end{pmatrix}, $$
$$S^i=\begin{pmatrix}2i+1 & -2i & i \\ 2i & -(2i-1) & i \\ 0 & 0 & 1
\end{pmatrix}, $$
$$S^iG= \begin{pmatrix}2i+1 & 1 & i \\ 2i & 1 & i \\ 0 & 0 & 1
\end{pmatrix},$$
$$(SG)^i=  \begin{pmatrix}a(i) & b(i) & c(i) \\ d(i) & a(i-1) & b(i) \\ 0 & 0 & 1
\end{pmatrix}$$
where 
\begin{itemize}
    \item $a(n)$ satisfies $a(n)=4a(n-1)-a(n-2)$ with $a(0)=1$ and $a(1)=3$ (it is, up to a shift, sequence A001835 in \cite{Sl}).
    \item $b(n)$ satisfies $b(n)=4b(n-1)-b(n-2)$ with $b(0)=0$ and $b(1)=1$ (it is sequence A001353 in \cite{Sl}).
    \item $c(n)$ satisfies $c(n)=5c(n-1)-5b(n-2)+c(n-3)$ with $c(1)=1$ and $c(k)=0$ for $k\leq 0$ (it is sequence A061278 in \cite{Sl}).
    \item $d(n)$ satisfies $d(n)=4d(n-1)-d(n-2)$ with $d(0)=0$ and $d(1)=2$ (it is sequence A052530 in \cite{Sl}).
\end{itemize}
\end{lem}
\proof
The assertions are easily provable by induction. 
\endproof

Now consider the subset of $\mathbb Z\times\mathbb Z$ whose points have equal positive coordinates
$$A=\{(k,k)\,|\,k\in\mathbb Z\,k\geq 0 \}.$$
Clearly $$A=\cup_{i\geq 0}\{S^i(1,1)\}\cup\{(0,0)\}.$$
Hence, $$F_{k,k}=1$$ if $k\geq0.$
Moreover, by the previous lemma, $$\cup_{i\geq 0} G^i(A)=\{(kd,d)\,|\,k,d\in\mathbb N\,,k>0\}.$$
This implies the following proposition. 

\begin{prop}\label{divisorn}
If $n\geq 1$ and $k$ divides $n$ then $F_{n,k}=1.$
\end{prop}
% \proof
% We prove the assertion by induction. The base step consists in the trivial verification that $F_{1,1}=1.$
% Suppose by induction that the assertion is true for every natural number less than $n.$
% Let $n=a\cdot k.$ If $a=1,$ we have $F_{n,k}=F_{k,k}=1.$
% Otherwise $a>1$ and, by Lemma \ref{lemma},
% $F_{n,k}=F_{ak,k}=F_{(a-1)k,h}$ where $h\equiv k \mod{2(a-1)k+1}$ and $0\leq h\leq 2(a-1)k. $ Since $2(a-1)k+1>k,$ $h=k$ and we have $F_{n,k}=F_{(a-1)k,k}.$
% This last term is equal to $1$ by the inductive hypothesis. 
% \endproof

The previous proposition implies that $$F_n\geq \tau(n),$$
where $\tau$ is the number-of-divisors function
(see e.g. \cite{apostol1998introduction}). We will see below that the previous inequality can be substantially improved.

\section{A modular recursion}

Now we prove a lemma that show how the entries of the matrix $F$ are related to each others through a modular recursion of the indices. 

\begin{lem}\label{lemma}
Consider the matrix $F.$ Then
\begin{equation}\label{ric}
F_{n,k}=F_{n-k,h},
\end{equation}
where $h$ is the remainder of $2n-2k+1$ in the division by $k,$ i.e. $$h\equiv k \mod{2n-2k+1}$$ and $0\leq h\leq 2n-2k.$    
Equation (\ref{ric}) together with the initial condition $F_{0,0}=1$  characterizes the matrix $F.$
\end{lem}
\proof
As recalled above, if $F_{n,k}=1$ then $1\leq k\leq n$ or $k=n=0.$
The lemma is clearly true for the entries of the form $F_{n,n},$
$n\geq 0,$ i.e. those obtained applying only the operation $S$ to $(1,1).$
In fact, in this case, $F_{n,n}=1,$ $h=0,$ $n-k=0$ and $F_{0,0}=1.$ 

Now choose a pair $(n,k)$ with $k<n.$ We want to show that 
$(n,k)=S^iG(n-k,h),$ $i\geq 0.$
By Lemma \ref{matrices}, $(n,k)=S^iG(x,y)$ if and only if
$n=(2i+1)x+y+i$ and $k=2i x+y+i.$ Hence $n-k=x,$ $2n-2k+1=2x+1$ and $h=y.$

Note that $F_{n,k}=1$ with $k<n$ if and only if $(n,k)=w(1,1)$ where $w\in H$ is a non-empty word with at least one letter equal to $G.$
Let $w=S^{i_l}GS^{i_{l-1}}G\ldots S^{i_1}GS^{i_0},$ where $l_j\geq 0$ for all $0\leq j\leq l.$ Set $$(x,y):=S^{i_{l-1}}G\ldots S^{i_1}GS^{i_0}(1,1).$$ Then $(n,k)=S^{i_l}G(x,y).$ Thus $F_{n,k}=F_{x,y}.$
\endproof

\section{The number-of-divisors function and the matrix $F$}

\begin{thm}\label{sdiv}
For every $n\geq 1$ and for every $x$ such that 
\begin{itemize}
    \item $0\leq x\leq n-1$ and
    \item there exists a divisor $h$ of $x$ with $2x+1\,|\,2n-2h+1,$
\end{itemize}
it holds $F_{n,n-x}=1.$
\end{thm}
\proof
If $x=0$ the proposition is trivial.
Let $x\geq 1$ be an integer. Then, by Lemma \ref{lemma}, we have
$F_{n,n-x}=F_{x,h}$ where $h\equiv n-x \mod{2x+1}.$ Now, by Proposition \ref{divisorn}, we have that $F_{x,h}=1$ if $h\,|\,x.$
Moreover $h\equiv n-x \mod{2x+1}$ if and only if $2x+1\,|\, n-x-h.$ Since $2x+1$ is odd this last condition is equivalent to $2x+1\,|\,2n-2x-2h$ and this in turn is equivalent to $2x+1\,|\, 2n-2h+1.$
Hence, if $h|x$ and $2x+1\,|\, 2n-2h+1,$ we have $F_{n,n-x}=1.$
\endproof

We denote by $a_n$ the number of entries $F_{n,n-x}$ of row $n$ of the matrix $F$ such that $0\leq x \leq n-1$ and there exists a divisor $h$ of $x$ with $2x+1|2n-2h+1.$
Clearly $$a_n \leq F_n,$$ for every $n.$

\begin{thm}\label{disug}
We have $$a_n\geq \max\{\tau(n),\,\tau(2n-1),\,\tau_o(n+1)\},$$
more precisely
$$a_n\geq \tau(n)+\tau(2n-1)+\tau_o(n+1)-3-\delta_{n\equiv 0 \mod 2}-\delta_{n \equiv -1 \mod 3}$$ for every $n\geq 1,$
where $\tau$ is the number-of-divisors function and $\tau_o$ is the number-of-odd-divisors function.
\end{thm}
\proof
Let $n-x$ be a divisor of $n.$ Then $n=(n-x)j,$ where $j$ is an integer. Then $2x+1=2n-2\frac{n}{j}+1.$ Since $\frac{n}{j}=n-x$ is a divisor of $n,$ it divides also $x,$ hence $F_{n,n-x}$ is one of those entries counted by $a_n.$
Hence $a_n\geq \tau(n).$ Notice that, in this case, $x=0$ or $n-x\leq \frac{n}{2}.$

From the previous theorem, taking $h=1,$ it follows, in particular, that $F_{n,n-x}=1$ if $2x+1$ is a divisor of $2n-1.$
Hence $a_n\geq \tau(2n-1)$ and $x=n-1$ or $2x+1\leq \frac{2n-1}{2}.$ In the last case, $x\leq \frac{n-1}{2}$ and 
$n-x\geq \frac{n+1}{2}.$ As a consequence
$a_n\geq \tau(n)+\tau(2n-1)-2,$ where the 2 in the right-hand side takes into account the fact that $F_{n,n}$ and $F_{n,1}$ have been counted two times. 

Similarly, taking $h=x$ in the previous theorem, it follows that $F_{n,n-x}=1$ when
$x$ is such that $2x+1\,|\, 2n-2x+1.$ This is equivalent to $2x+1\,|\,2n+2$ which, since $2x+1$ is odd, is in turn equivalent to $2x+1\,|\, n+1.$
Hence $F_{n,n-x}=1$ if $2x+1$ is an odd divisor of $n+1$ and
$a_n\geq \tau_o(n+1).$
Here $x=\frac{n}{2}$ or $2x+1\leq \frac{n+1}{2}.$   In the last case
$x\leq \frac{n-1}{4}$ and $n-x\geq \frac{3n+1}{4}.$

Notice that if $2x+1$ divides $n+1$ and $2n-1$ it divides also $2(n+1)-(2n-1)=3$ hence $x=0$ or $x=1.$ If $x=1,$ $2x+1=3$ divides $n+1$ if and only if it divides also $2n-1.$

Hence $$a_n\geq \tau(n)+\tau(2n-1)+\tau_o(n+1)-3-
\delta_{n\equiv 0 \mod 2}-\delta_{n \equiv -1 \mod 3}$$
where, in the right hand side, the 3 takes into account the fact that $F_{n,n}$ has been counted three times and $F_{n,1}$ has been counted two times, the $\delta_{n\equiv 0 \mod 2}$ takes into account the fact that $F_{n,\frac{n}{2}}$ has been counted two times if $n$ is even and $\delta_{n \equiv -1 \mod 3}$ takes into account the fact that $F_{n,n-1}$ has been counted two times if $n+1$ is divisible by 3.
\endproof

Since a number has 1 as its only odd divisor if and only if it is a power of 2, the previous theorem shows that Conjecture \ref{conj1} is proved for every $n> 7$ except those primes $p$ of the form $p=2^q-1$ such that $2p-1$ is also prime. Notice that the fact that $F_n\geq 3$ for every $n\neq 2^q-1$ follows also from Remark 1 in \cite{test1}.

A prime $p$ of the form $2^q-1$ is said to be a \textit{Mersenne prime} (see e.g. \cite{apostol1998introduction}). 
The Mersenne primes $p$ such that $2p-1$ is also prime appears in \cite{Sl} in sequence A167917.
It is well known that, if a prime $p$ has the form $2^q-1,$ then also $q$ is prime.

\begin{thm}\label{diop}
The number $a_n$ is one plus the number of solutions of the Diophantine equation 
\begin{equation}\label{equa}
n=2xyz +yz +x+ y, \quad \mbox{ with }\,y,z\geq 1\,\mbox{ and }\,x \geq 0.
\end{equation}
Moreover, $$a_n=1+\sum_{0\leq j<n}D_{2j+1}(n-j),$$
where $D_m(n)$ is the number of divisors $d>1$ of $n$ congruent to 1 $\mod m.$
\end{thm}
\proof
Let $n$ be fixed. The number of elements of the form $F_{n,n-x}$
with $x\neq 0,$ such that there exists an $h$ with $h|x$ and $2x+1|2n-2h+1$
is equal to the number of solutions $h,k,j\geq 1$ to the equation
\begin{equation}
2n-2h+1=j(2hk+1)    
\end{equation}
i.e. $$2n=2jhk+2h+j-1.$$ In this equation $j$ must be odd and hence $j=2\hat j +1.$ So we get the equation
\begin{equation}\label{equa2}
n=2\hat j hk +hk +h+ \hat j, \quad \mbox{ with }\,h,k\geq 1\,\mbox{ and }\,\hat j \geq 0.
\end{equation}
Hence $a_n$ is equal to the number of solutions to this equation increased by one since we have to take into account the case with $x=0.$

Consider now Equation \ref{equa2}. It is equivalent to $$n-\hat j=h((2\hat j  +1)k+1).$$ Hence,
if $n$ and $\hat j\geq 0$ are fixed, there is a correspondence between the solutions $h,k\geq 1$ to the last equation and the number of divisors $d=(2\hat j  +1)k+1>1$ of $n-\hat j$ congruent to one $\mod 2\hat j +1.$ 
\endproof

Now we prove Conjecture \ref{conj1}.

\begin{thm}\label{Thmconj1}
$F_n\geq 3$ for all $n>7.$
\end{thm}
\proof
As recalled above, the assertion follows from the previous results for every $n> 7$ except those primes $p$ of the form $p=2^q-1,$ $q$ a prime, such that $2p-1$ is also prime.

Let $p$ be a number with these properties such that $F_p=2.$
Since $F_n\geq a_n$ for every $n$ and since $a_n\geq 2$ by Theorem \ref{disug}, we have $a_p=2.$
Moreover $$a_p\geq 1+D_1(p)+D_3(p-1).$$ Clearly $D_1(p)=1.$
Hence $D_3(p-1)=0.$ We want to show that it is impossible. 
The integer $p-1$ has no divisors $d>1$ such that
$d\equiv 1 \mod{3}$ if and only if $p-1=3^r p_1^i$ where $p_1$ is a prime with $p_1\equiv -1 \mod{3},$  $i=0$ or $i=1$ and $r\geq 0.$ In fact, if between the prime factors of $p-1$ there were more than one congruent to $-1$ $\mod{3}$ or at least one congruent to $1$ $\mod{3},$ then $p-1$ would have at least one divisor congruent to $1$ $\mod{3}.$

Since $p=2^q-1,$ we have $$2(2^{q-1}-1)=3^r p_1^i$$
which implies $p_1=2,$ $i=1$ and $2^{q-1}-3^r=1.$
The only solutions to the previous equation in positive integers $q$ and $r$ are $q=3$ and $r=1.$
In fact, in $2002,$ Mih\u{a}ilescu proved that the only solution to the Diophantine equation $x^a-y^b=1,$ with $x,y,a,b>1$  is $x=3,$,$y=2,$ $a=2,$ $b=3,$ thus solving the celebrated Catalan's conjecture (see \cite{schoof2010catalan},
 the solution of the particular case of this conjecture with $x=2$ and $y=3$ is attributed to Gersonides).

If $q=3$ and $r=1,$ we get $p=7,$ 
whereas we are considering a number $p>7.$
This concludes the proof.

\endproof

\section{Other properties of the matrix $F$}

In this section we investigate further properties of matrix $F.$

\begin{thm}
 Matrix $F$ has periodic diagonals. In particular, the $a$-th  subdiagonal has period $2a+1.$
 
 Matrix $F$ has periodic columns. In particular, column $a$ has period $a$.
\end{thm}
\proof
To prove the first part of the theorem, fix $a\in \mathbb N$ and consider the elements $F_{n,n-a},$ $n> a,$ of the matrix $F.$
These elements constitute the $a$-th subdiagonal. By Lemma \ref{lemma} we have $F_{n,n-a}=1$ if and only if $F_{a,h}$
where $h\equiv n-a \mod{2a+1}.$
Since $a$ is fixed, this proves that the $a$-th subdiagonal is periodic with period $2a+1.$

To prove the second part, fix $a\in \mathbb N$ and consider 
the elements $F_{n,a},$ $n\geq a.$ These elements constitutes the $a$-th column of $F.$
By Lemma \ref{lemma}  $F_{n,a}=1$ if and only if $F_{n-a,h}$
where $h\equiv a \mod{2n-2a+1}.$ If $n$ is sufficiently large, this implies $h=a$ and hence the $a$-th column has period $a.$ 
\endproof

\begin{thm}\label{kdj}
For each quadruple $(k,d,i,j)\in \mathbb N^4,\, k>0,$ we have 
$$F_{kd+id+j(i+1)(2(k-1)d+1),d+j(2(k-1)d+1)}=1.$$
\end{thm}
\proof
Consider the set $A:=\{(kd,d)\,|\,k,d\in \mathbb N\,,k>0\}.$
We have $F_{x,y}=1$ for all $(x,y)\in A$ by Proposition \ref{divisorn}.
By Lemma \ref{lemma}, we have
$$\cup_{i,j\geq 0}  G^i  S^j(A)=$$ $$\{kd+id+j(i+1)(2(k-1)d+1),d+j(2(k-1)d+1))\,|\,k,d,j,l\in\mathbb N\,,k>0\}.$$
This concludes the proof. 
\endproof

The previous theorem implies the following Corollary.

\begin{cor}\label{tredue}
For each $t \in \mathbb N,$ $F_{3t+2,2t+2}=1$ and $F_{5t+4,2t+2}=1.$
\end{cor}
\proof
For the first part, substitute $i=0,d=1,j=1$ and $k-1=t$ in the previous theorem.
% For the first part, by Lemma \ref{lemma}, we get $F_{3j+2,2j+2}=1$ if and only if $F_{j,h}=1$
% where $h\equiv 2j+2 \mod{2j+1}$ i.e. $h=1.$ Since $F_{j,1}=1,$
% we get the assertion. 
For the second part, by Lemma \ref{lemma}, we get $F_{5j+4,2j+2}=1$ if and only if $F_{3j+2,h}=1$
where $h\equiv 2j+2 \mod{6j+4}$ i.e. $h=2j+2.$ Since $F_{3j+2,2j+2}=1$ by the first part,
we get the assertion. 
\endproof

%The following theorem generalizes the first assertion of the %last one.

% \begin{thm}\label{kdj}
% For each triple $(k,d,j)$ of positive integer,
% $$F_{kd+j(2(k-1)d+1),d+j(2(k-1)d+1)}=1$$
% \end{thm}
% \proof
% We use again Lemma \ref{lemma}.
% $$F_{kd+j(2(k-1)d+1),d+j(2(k-1)d+1)}=1$$
% if and only if
% $$F_{(k-1)d,h} =1$$
% with $h\equiv d+j(2(k-1)d+1) \mod{2(k-1)d+1},$ i.e. $h=d.$
% We get the assertion since $F_{(k-1)d,d}=1$ by Theorem \ref{divisorn}.
% \endproof

\section{Conjectures about the matrix $F$}

In this section we formulate others conjectures about $F$ and explain their relation with Conjecture \ref{conj1}. To this aim we need to introduce the notion of track vector.

Denote by $\phi$ the map that associates the pair of integers $(n,k)$ the pair $(n-k,h)$ where $h\equiv k \mod 2n-2k+1$ and $0\leq h\leq 2n-2k.$ It follows from the proof of Lemma \ref{lemma} that $\phi(n,k)=(n-k,h)$ if and only if there exists an $i\in \mathbb N$ such that $S^iG(n-k,h)=(n,k).$
As a consequence, $(n,k)=S^{i_l}GS^{i_{l-1}}G\ldots S^{i_1}G(m,m),$ where $i_1,\ldots,i_l,m\in \mathbb N,$ if and only if $\phi^{i_l}(n,k)=(m,m).$ 

% \begin{lem}
% Let $(n,k)$ be a pair of different integers with $F_{n,k}=1.$
% The number of times one need to apply the operation described in Theorem \ref{sdiv} to get a pair of equal integers $(m,m)$ is equal to the number of times it is necessary to apply an operation  of the form $\hat S^i\hat G,$ with $i\geq 0,$ to  $(m,m)$ to get $(n,k).$
% \end{lem}
% \proof
% If the element $(n,k,1)\in \mathbb Z^3$ with $n\neq k$ and $F_{n,k}=1$, it has been obtained by $(1,1,1)$ multiplying by the matrix
% $$S^{i_l}GS^{i_{l-1}}G\ldots S^{i_1}GS^{m}$$ with $i_J\geq 0$ and $l\geq 1$ (if $l=0,$ $n$ would be equal to $k.$)
% Notice that $S^{m}(1,1,1)=(m,m,1).$
% By the previous lemma, it follow that $$S^iG= \begin{pmatrix}2i+1 & 1 & i \\ 2i & 1 & i \\ 0 & 0 & 1
% \end{pmatrix},$$ for all $i\geq 0.$
% Hence $$S^iG(x,y,1)=((2i+1)x+y+i,2ix+y+i,1).$$
% Applying to this last vector Theorem \ref{sdiv} we get that 
% $$F_{(2i+1)x+y+i,2ix+y+i}=1 $$ if and only if $F_{x,y}=1.$
% Hence an application of Theorem \ref{sdiv} correspond to a multiplication by the matrix $$(S^iG)^{-1}.$$
% \endproof

The number of operations of the form $S^iG$ needed to get a given element is related to the breadth of an element. Here we recall the definition of breadth and of track vector. 
Following \cite{test1}, the\textit{track vector} of an element $(n,k)$ is defined as the vector $(i_0+1,\ldots,i_l+1)$ where $$(n,k)=S^{i_l}GS^{i_{l-1}}G\ldots S^{i_1}GS^{i_0}(1,1).$$
In this case, the \textit{breadth} of $(n,k)$ is equal to $l.$
Since $S^{i_0}(1,1)=(i_0,i_0),$ the breadth of $(n,k)$ is equal to the number of times it is necessary to apply the map $\phi$ to $(n,k)$ to get an entry of the form $(m,m).$

\begin{thm}
The elements of $F$ appearing in Theorem \ref{sdiv} with $x\neq 0,$
i.e. the elements $F_{n,n-x},$ $x\neq 0,$ such that there exists an $h$ with $h|x$ and $2x+1|2n-2x+1$ are precisely the elements in row $n$ with track vector of the form $(a,(1)^p,b),$ with $a,b,p\geq 1.$

Moreover, if $n$ is fixed, the number of elements $F_{n,n-x},$ $x\neq 0,$ such that there exists an $h$ with $h|x$ and $2x+1|2n-2x+1$ 
is equal to the number of elements in row  $n$ with track vector $((1)^p,a,(1)^q),$ where $p,q\geq 0,$ $a\geq 2.$
\end{thm}
\proof
By the previous observations and by the proof of Theorem \ref{sdiv},
we have that the elements  $F_{n,n-x},$ $x\neq 0,$ such that there exists an $h$ with $h|x$ and $2x+1|2n-2x+1$ are precisely the elements of the form 
$$S^{b-1}G^{p+1}S^{a-1}(1,1)$$ with $a,b,p\geq 1$ i.e. those with track vector $(a,(1)^p,b).$

As in the previous Sections, we denote by $a_n$ the number of elements in row $n$ of the form described in the proposition. 
It follows from Theorem \ref{diop} that $a_n-1$ is equal to the number of solutions to Equation \ref{equa}.
By the first part of Lemma 8 in  (\cite{test1}),  
the number of elements of the form $(n,r)$ with track vector $((1)^p,a,(1)^q)$
is equal to $$ap(2q+2)+a(q+1)-p(2q+1)=2ap(q+1)+a(q+1)-2p(q+1)+p.$$
If $n$ is fixed and we substitute $q+1=\hat q$ the number of such elements is equal to the number of solutions of the equation $$n=2ap\hat q+a\hat q-2p\hat q+p=2p\hat q (a-1)+(a-1)\hat q +\hat q +p.$$
Set $\hat a =a-1.$ We get the equation $n=2p\hat q \hat a +\hat a \hat q+\hat q +p,$ where $\hat a , \hat q\geq 1,$ and $p \geq 0.$
This equation coincides with  Equation \ref{equa}. 
\endproof

\begin{cor}
The number of elements in row $n$ with track vector of the form $(a,(1)^p,b),$ with $a,b,p\geq 1$ is equal to the number of elements in row $n$ with track vector of the form $((1)^p,a,(1)^q),$ where $p,q\geq 0,$ $a\geq 2.$
Moreover this common value is $a_n-1.$
\end{cor}

We conjecture that the sequence $a_n$ tends to infinity.

\begin{conj}\label{aninf}
$$a_n\to \infty,$$ more precisely $a_n\geq \lfloor \log(n)\rfloor-1.$
\end{conj}

The inequality of the conjecture originates from numerical evidences. Notice that it is not true that $a_n\geq \lfloor \log(n)\rfloor.$ In fact $a_{18007}=8$ but $\lfloor \log(18007)\rfloor=9.$

Notice that the previous conjecture implies Conjecture \ref{conj1} and Theorem \ref{Thmconj1}.

Moreover, by the paper \cite{test1}, it follows that the set of elements in row $n$ with track vector $((1)^p,a,(1)^q),$ where $p,q\geq 0,$ $a\geq 2$ corresponds bijectively with a subset of the  \textit{elementary partitions}, which in turn are a subset of the set of partitions with $n$ subpartitions. Hence, if $s_n$ is the number of partitions with $n$ subparitions, the previous conjecture implies also that $s_n \to \infty.$
Notice that sequence $\{s_n\}_{n\in\mathbb N}$ is sequence A116473 in \cite{Sl}, where it is reported that it is conjectures that $s_n\to \infty$. Thus Conjecture \ref{aninf} would imply also this conjecture present in \cite{Sl}.

Another conjecture suggested by strong numerical evidences is the following.
\begin{conj}
Let $c_n$ be the number of elements in row $n$ with breadth $3.$ Then $c_n\to \infty.$
\end{conj}

We conclude this Section with a conjecture about the possible positions of elements created at Step $i$ inside the matrix $F.$

\begin{conj}
Consider the elements of matrix $F$ that are created at Step $i.$
Let $r_i$ be the maximal index of a row of matrix $F$ containing such an element. Then $\{r_{2j+1}\}_{j\in \mathbb N}$ is sequence A061278 in \cite{Sl} and
$\{r_{2j}\}_{j\in \mathbb N}$ is sequence A001571 in \cite{Sl}.
Moreover, if $i$ is even, row $r_i$ of $F$ contains two elements created at Step $i.$ On the other hand, if $i$ is odd, row $r_i$ contains one element created at Step $i.$
\end{conj}
%\proof
%Induction.
%\endproof

\begin{example}
Let $i=6.$ The maximal row containing an element created at Step 6 is row 35. Notice that 35 is the third element (avoiding the first zero) of sequence A001571 in \cite{Sl}.
Moreover, row 35 of $F$ contains two elements created at Step 6, precisely $F_{35,15}$ and $F_{35,26}.$
\end{example}

\section{Upper bound for $F_n$}

We conclude the paper improving the best known upper bound for $F_n.$

Corollary 5 in \cite{BaBoCaCo} states that $F_n\leq \min\{n,\phi(2n+1)\}$ where $\phi$ is the Euler totient function, see e.g. \cite{apostol1998introduction}.
In fact $F_n\leq n$ since the matrix $F$ is lower triangular. Moreover, it is shown in \cite{BaBoCaCo}, that if $F_{n,k}=1$ then $\gcd(k,2n+1)=1$ (this can also be shown easily using the recursive definition of $F$ used in this paper). Hence $F_n\leq \phi(2n+1).$

It is possible to slightly improve this bound in the following way. Consider the same notation of the proof of Theorem \ref{free}. Since $G(1,1)=(2,1)\in X_1,$ $S(1,1)=(2,2)\in X_2,$
$G(X_1\cup X_2)\subseteq X_1$ and $S(X_1\cup X_2)\subseteq X_2$
we have that every element $(n,k),$ $n>1,$ such that $F_{n,k}=1$ is contained in $X_1\cup X_2.$  In particular, $k<\frac{n}{2}+\frac{1}{4}$ or $k>\frac{2n}{3}+\frac{1}{3}.$
Hence $$F_n\leq n-\left(\frac{2n}{3}+\frac{1}{3}-\left(\frac{n}{2}+\frac{1}{4} \right)-1\right)=\frac{5n}{6}-\frac{13}{12}.$$
Thus we get $$F_n\leq \min\left\{\frac{5n}{6}-\frac{13}{12},\phi(2n+1)\right\}.$$

% \begin{conj}
% $$\lfloor \sqrt{n}\rfloor -1 \leq F_n\leq \frac{n}{log(log(n))+\frac{3}{log(log(n))}} $$
% for every $n.$
% \end{conj}
% The right hand side of the previous inequality is a well-known lower bound for the Euler totient function. 

\addcontentsline{toc}{section}{Bibliography}
\bibliographystyle{plain}
\bibliography{BIBLIOGRAFIA}

\end{document}